\documentclass[12pt]{article}
\usepackage{graphicx}
\usepackage{amsfonts}
\usepackage{amssymb}
\usepackage{geometry}

\input{tcilatex}
\begin{document}

\begin{center}
{\huge The Graph of the Hypersimplex}

\bigskip

\bigskip

\bigskip

{\Large Fred J. Rispoli}

\bigskip

\textit{Department of Mathematics and Computer Science}

\textit{Dowling College, Oakdale, NY 11769}

\bigskip

\bigskip

\bigskip

\bigskip

\bigskip

\bigskip ABSTRACT
\end{center}

The $(k,d)$-hypersimplex is a $(d-1)$-dimensional polytope whose vertices
are the $(0,1)$-vectors that sum to $k$. When $k=1$, we get a simplex whose
graph is the complete graph $K_{d}$. Here we show how many of the well known
graph parameters and attributes of $K_{d}$ extend to a more general case. In
particular we obtain explicit formulas in terms of $d$ and $k$ for the
number of vertices, vertex degree, number of edges and the diameter. We show
that the graphs are vertex transitive, hamilton connected, obtain the clique
number and show how the graphs can be decomposed into self-similar
subgraphs. The paper concludes with a discussion of the edge expansion rate
of the graph of $\ $a $(k,d)$-hypersimplex which we show is at least $d/2$,
and how this graph can be used to generate a random subset of $%
\{1,2,3,...,d\}$ with $k$ elements.

\bigskip

\bigskip

\bigskip

\section{Introduction}

The $(k,d)$-hypersimplex, denoted by $\Delta _{d,k}$, is defined as the
convex hull of all ($0,1)$-vectors in $%
%TCIMACRO{\U{211d} }%
%BeginExpansion
\mathbb{R}
%EndExpansion
^{d}$ whose nonzero elements sum to $k$. Hypersimplices are $(d-1)$%
-dimensional polytopes that appear in various algebraic and geometric
contexts (e.g., see [6]). The polytope $\Delta _{d,k}$ can also be defined
as a "slice" of the $(d-1)$-hypercube located between the two hyperplanes $%
\sum x_{i}$ $=d-1$ and $\sum x_{i}$ $=d$ in $%
%TCIMACRO{\U{211d} }%
%BeginExpansion
\mathbb{R}
%EndExpansion
^{d}$. A classical result implied by the work of Laplace [7], is that the
normalized volume of this polytope equals the Eulerian number $A_{k,d-1\text{%
.}}$\ Hypersimplices have also appeared in the theory of characteristic
classes and\ Gr\"{o}bner bases (for more details on this and on polytopes in
general, see [3] and [12].) The graph of the hypersimplex $\Delta _{d,k}$,
denoted by $G_{d,k}$, is the graph consisting of the vertices and edges of $%
\Delta _{d,k}$. This graph is also known as the Johnson graph and $G_{d,k}$
provides an example of a family of "distance-regular" graphs of unbounded
diameter, which are also a special type of Coxeter graph [1].

For the case $k=1$, the graph $G_{d,1}$ is the complete graph $K_{d}$ whose
role is fundamental in Graph Theory.  Compared to the complete graphs, the
properties of the closely related hypersimplex graphs are not very well
known. Here we show how many of the parameters of $K_{d}$ are extended to $%
G_{d,k}.$ For example, $K_{d}$ has $d$ vertices, is regular of degree $d-1$
and has diameter $1$. We show that $G_{d,k}$ has $\left( 
\begin{array}{c}
d \\ 
k%
\end{array}%
\right) $ vertices, is regular of degree $k(d-k)$, and has diameter $k$, for 
$k$ $\leq \frac{d}{2}$. We also characterize adjacency, show that $G_{d,k}$
is vertex transitive, obtain an explicit formula for the number of edges,
determine the clique number for $G_{d,k}$, and study varous connectivity
properties. In addition, since the number of vertices in $G_{d,k}$ is $%
\left( 
\begin{array}{c}
d \\ 
k%
\end{array}%
\right) $, it is natural to ask how $G_{d,k}$ can be decomposed into
subgraphs whose vertex counts satisfy Pascal's Identity $\left( 
\begin{array}{c}
d \\ 
k%
\end{array}%
\right) =$ $\left( 
\begin{array}{c}
d-1 \\ 
k-1%
\end{array}%
\right) +\left( 
\begin{array}{c}
d-1 \\ 
k%
\end{array}%
\right) $. We show that this leads to a recursive decomposition of\ $G_{d,k}$%
\ into self similar subgraphs. The paper concludes with a discussion of edge
expansion properties of $G_{d,k}$ and random walks on $G_{d,k}$ that may be
used to generate random subsets of $\{1,2,3,...,d\}$ of size $k$.

\bigskip

\bigskip 

\section{\protect\bigskip The vertices and edges of $G_{d,k}$}

\bigskip

A polytope contained in $%
%TCIMACRO{\U{211d} }%
%BeginExpansion
\mathbb{R}
%EndExpansion
^{d}$ is called $(0,1)$-\textit{valued} if all of its vertices are vectors
having coordinates that are all either $0$ or $1$. Given any convex polytope 
$P$, two distinct vertices $x\neq y$ in $P$ are \textit{adjacent} if for
every $\lambda $\ satisfying $0<\lambda \ <1$, it holds that\ $\lambda
x+(1-\lambda )y$ can not be expressed as a convex combination of other
vertices in $P$. A graph is said to be \textit{regular} of degree $r$ if
every vertex in the graph has degree $r$. Let $x$ and $y$ be points in $%
%TCIMACRO{\U{211d} }%
%BeginExpansion
\mathbb{R}
%EndExpansion
^{d}$, then $x\cdot y$ is the inner product $\dsum%
\limits_{i=1}^{d}x_{i}y_{i} $.

\bigskip

\bigskip

\textbf{Proposition 1} \textit{For} $1\leq k<d$ and $d\geq 4$:

(a) \textit{The number of vertices in }$G_{d,k}$ \textit{is} $\left( 
\begin{array}{c}
d \\ 
k%
\end{array}%
\right) .$

(b) \textit{Two distinct vertices x and y of }$G_{d,k}$ \textit{are adjacent
if and only if }$x\cdot y=k-1$.\bigskip

\bigskip

\textbf{Proof.} (a) The count follows from the fact that there is an obvious
one-to-one correspondence between the subsets of $\{1,2,...,d\}$ with $k$
elements and the number of $0,1$ $d$-vectors with exactly $k$ ones.

(b) If $k=1$, then $G_{d,k}$ is the complete graph $K_{d}$ and the result
holds since all vertices in $K_{d}$\ are adjacent. So assume that $k\geq 2$,
and suppose that $x\cdot y<k-1$. Then there exists $p,$ $q,$ $r,$ $s\in
\{1,2,...,d\}$ such that $x_{p}=1$, $x_{q}=1$, $x_{r}=0$, $x_{s}=0$, and $%
y_{p}=0$, $y_{q}=0$, $y_{r}=1$, $y_{s}=1$. Define vertices $u$ and $v$ as
follows: $x_{i}=u_{i}$, for all $i$ except $i=q,$ $r$, which satisfy $%
u_{q}=0 $ and $u_{r}=1$, and $y_{i}=v_{i}$, for all $i$ except $i=q,$ $r$,
which satisfy $v_{q}=1$ and $v_{r}=0$. Then $\frac{1}{2}x+\frac{1}{2}y$ $=%
\frac{1}{2}u+\frac{1}{2}v,$ so $x$ and $y$ are not adjacent.

Next observe that $x\cdot y>k$ is impossible, and $x\cdot y$ $=k$ implies
that $x=y$, which is a contradiction. So suppose that $x\cdot y$ $=k-1$ and
that there exists $z^{1}$, $z^{2}$, ... , $z^{n}$ such that $\lambda
x+(1-\lambda )y=\dsum\limits_{j=1}^{n}\alpha _{j}z^{j}$. Since $x$ and $y$
both have exactly $k$ ones, $d-k$ zeroes, and $x\cdot y$ $=k-1$, $x$ and $y$
must be equal for all indices except for $2$. This implies that the $z^{j}$
are equal on all indices except for $2$, and hence $n=2$. Consequently, we
must have $\{x,y\}=\{z^{1},z^{2}\}$. \ $\blacksquare $

\bigskip

\bigskip

\textbf{Proposition 2} (a) \textit{The graph }$G_{d,k}$ \textit{is regular
of degree }$k(d-k)$.

\bigskip (b) \textit{The number of edges in} $G_{d,k}$ \textit{is} $\frac{d!%
}{2(k-1)!(d-k-1)!}.$

\bigskip

\textbf{Proof.} (a) Let $x$ be any vertex in $G_{d,k}$. By Proposition 1, a
vertex $y$ is adjacent to $x$ if and only if $x\cdot y=k-1$. So $y$ must
have $k-1$ of the $k$ ones in $x$. There are $k$ ways for this to happen. In
addition, one of the $d-k$ indices for which $x_{i}=0$ must be a one for $y$%
. There are $d-k$ ways for this to happen. Hence the number of vertices
adjacent to $x$ is $k(d-k)$.

(b) The count follows from the well known Handshaking Lemma [11] and the
fact that the sum of the vertex degrees in $G_{d,k}$\ is given by $\left( 
\begin{array}{c}
d \\ 
k%
\end{array}%
\right) k(d-k)$. \ \ $\blacksquare $

\bigskip

\bigskip

A graph $G$ with vertices $V(G)$ is called \textit{vertex transitive} if
given any two vertices $x$ and $y$ there is an automorphism $%
f:V(G)\rightarrow V(G)$ such that $f(x)=y$. It is known that graphs of
Platonic solids and Archimedean solids are vertex transitive, as well as $%
K_{d}$ and the complete bipartite graph $K_{d,d}$. We now show that this
property is also true for $G_{d,k}$.

\bigskip

\bigskip \textbf{Proposition 3} \textit{For} $1\leq k<d,$ \textit{the graph }%
$G_{d,k}$ \textit{is vertex transitive}.

\bigskip

\textbf{Proof.} Given vertices $x$ $=(x_{1}x_{2}...x_{d})$ and $%
y=(y_{1}y_{2}...y_{d})$ define $f:V(G_{d,k})\rightarrow V(G_{d,k})$ as
follows. (An example illustrating the construction of $f$ is given
immediately after the proof.) If $x_{i}=y_{i}$, $f$ takes $i$ to $i$; i.e., $%
f$ takes the $i$'th digit in $x$ to the $i$'th digit in $f(x)$. Now consider
the set $D(x,y)=\{i:x_{i}\neq y_{i}\}$. Let $p$ be the number of elements in 
$D(x,y)$ for which $x_{i}=0$, and let $q$ be the number for which $x_{i}=1.$
Then the number of elements in $D(x,y)$ for which $y_{i}=1$ must be $p$ and
the number for which $y_{i}=0$ is $q$. But $\sum x_{i}$ $=\sum y_{i}$, so we
must have $p=q,$ and hence, $D(x,y)$ contains $2p$ indices. Let $i_{1}$\ be
the smallest index for $x$ in $D(x,y)$ with $x_{i_{1}}=1$ and $i_{2}$ the
smallest index in $D(x,y)$ with $x_{i_{2}}$ $=0$. Let $j_{1}$\ be the
smallest index for $y$ in $D(x,y)$ with $y_{j_{1}}=0$ and $j_{2}$ the
smallest index with $y_{j_{_{2}}}=$ $1$. Then we define $f$ such that the $%
i_{1}$th digit in $x$ becomes the $j_{2}$th digit in $f(x)$, and the $i_{2}$%
th digit in $x$ becomes the $j_{1}$th digit in $f(x)$. Since $D(x,y)$ has an
even number of indices, we may repeat this step as often as necessary.

Observe that $f(x)=y$. In addition, $f$ is its own inverse. Hence $f(x)=y$
implies that $f(f(x))=f(y)$, or simply $x=f(y)$. This gives $x\cdot
y=f(y)\cdot f(x)$, so $f$ is an adjacency preserving automorphism. $\ \
\blacksquare $

\bigskip

To illustrate the construction of $f$ suppose that $x=(110101101000)$ and $%
y\ =(100110010110)$. Then $f$ is given by $%
f(x_{1}x_{2}...x_{12})=(x_{1}x_{5}x_{3}x_{4}x_{2}x_{8}x_{10}x_{6}x_{11}x_{7}x_{9}x_{12}) 
$ $=y$. Furthermore, if $z=(101010101010),$ then $f(z)=(111000001110)$.

\bigskip

\section{Connectivity Properties}

\bigskip

The \textit{distance} between any two vertices $x$ and $y$, in a graph $G$
is the number of edges in a shortest path joining $x$ to $y$. The \textit{%
diameter} of $G$, $\delta (G)$, is the maximum distance amongst all pair of
vertices in $G$. A graph is called \textit{distance-regular} if it is a
regular graph such that, given any two vertices $x$ and $y$ at any distance $%
i\leq \delta (G)$, the number of vertices adjacent to $y$ and at a distance $%
j$ from $x$ depends only on $i$ and $j$, and not on the particular vertices.
A graph $G$ is called $d$-\textit{connected} if for every pair of vertices $%
x $ and $y$ there exists $d$ disjoint paths joining $x$ to $y$. A graph is
called \textit{hamilton connected} if every pair of distinct vertices is
joined by a path that passes through every vertex of $G$ exactly once. A
subset of vertices $H$\ is called a \textit{clique} in $G$ if there is an
edge in $G$ between every pair of vertices in $H$. The cardinality of the
largest clique in G is called the \textit{clique number} of $G$. It is easy
to show that $G_{d,k}$\ is isomorphic to $G_{d,d-k}$\ so in the following
propositions we restrict $k$ such that $1\leq k\leq \frac{d}{2}.$

\bigskip

\textbf{Proposition 4} \textit{Let }$1\leq k\leq \frac{d}{2}.$

\bigskip

(a) \textit{Given any two vertices} $x$ \textit{and} $y$ \textit{in} $%
G_{d,k} $, \textit{the distance between} $x$ \textit{and} $y$ \textit{is} $%
k-x\cdot y $.

(b) \textit{The diameter of} $G_{d,k}$\ \textit{is} $k$.

(c) $G_{d,k}$\ \textit{is a distance-regular graph}.

\bigskip

\textbf{Proof}. (a) Let $x\neq y$ be given. If $x$ and $y$ are adjacent,
then by Proposition 1, the distance between $x$ and $y$ is $1=k-x\cdot y$.
So assume that $x$ and $y$ are not adjacent and hence, $x\cdot y<k-1$. Since 
$x$ any $y$ both have $k$ ones and $d-k$ zeros, there exist indices $p$ and $%
q$ such that $x_{p}=1$, $y_{p}=0$, $x_{q}=0$ and $y_{q}=1$. Define the
vertex $z$ by $z_{i}=x_{i}$, for all $i$ except $p$ and $q$, where $z_{p}=0$%
, and $z_{q}=1$. Then $z$ has exactly $k$ ones and $x\cdot z=k-1$, so $x$
and $z$ are adjacent. Moreover, $x\cdot z=x\cdot y+1.$ We can repeat this as
often as necessary each time getting one step closer to $y$.

(b) Since $k\leq \frac{d}{2}$ there exists vertices $x$ and $y$ such that\ $%
x\cdot y=0$. By (a)\ this implies that the distance between $x$ and $y$ is $%
k $. Hence, the diameter must be $k$.

(c) Let $x$ and $y$ be vertices of $G_{d,k}$ whose distance is $i$. Then $x$
and $y$ have $k-i$ ones in common. Moreover, any vertex $z$ adjacent to $y$
at a distance $j$ from $x$, satisfies $y$ and $z$ have $k-1$ ones in common,
and $x$ and $z$ have $k-j$ ones in common. The number of vertices $z$
satisfying this depends only on $i$ and $j$. \ \ $\blacksquare $

\bigskip

\bigskip

\bigskip \textbf{Proposition 5} \ \textit{For} $2\leq k\leq \frac{d}{2}$:

(a) $G_{d,k}$ \textit{contains the complete graph} $K_{d-k+1}$ \textit{as a
subgraph}.

(b) \textit{The clique number of }$G_{d,k}$ \textit{is} $d-k+1.$

\bigskip

\bigskip

\textbf{Proof.} (a) Let $H$ be the subset of vertices whose first $k-1$
coordinates are all one. Then $H$ contains $d-k+1$ vertices, and every pair
of vertices $x$, $y$ in $H$\ satisfy $x\cdot y$ $=k-1$. Hence the subgraph
induced by $H$ must be $K_{d-k+1}$.

(b) Suppose, to obtain a contradiction, that the clique number of $G_{d,k}$\
is $w$ and $w>d-k+1$. Then there exists a subgraph isomorphic to $K_{p}$
where $p=d-k+2$. Let $x^{1},x^{2},...,x^{p}$ be the vertices of the subgraph.

If $x^{1},x^{2},...,x^{p}$\ all have $k-1$ ones in common, then without loss
of generality, we may assume that $x^{1},x^{2},...,x^{p}$ all have their
first $k-1$ digits equal to 1. Moreover, the last $d-(k-1)$ digits for\ each
of $x^{1},x^{2},...,x^{p}$ must all consist of zeros and exactly one 1.
Since the $x^{j}$ are all distinct, there are only $d-k+1$ possibilities for
this, implying that $p<d-k+2$, a contradiction.

Now suppose that $x^{1},x^{2},...,x^{p}$\ do not all have $k-1$ ones in
common, and that the first $k$ digits of $x^{1}$ are one. Observe that $k-1$
of the first $k$ digits of $x^{2},...,x^{p}$\ must be one since these
vertices must be adjacent to $x^{1}$. We show that no two of these vertices
have the same first $k$ digits. Suppose $x^{2}$ and $x^{3}$ have the same
first $k$ digits as illustrated below. Then, since $x^{1},x^{2},...,x^{p}$\
do not all have $k-1$ ones in common, there exists an $x^{4}$ whose first $k$
digits are different from those of $x^{2}$, also illustrated below. Notice
that $x_{k+1}^{4}$ $=1$, since $x^{4}$ and $x^{2}$ must be adjacent. But now 
$x^{4}$ and $x^{3}$ are not adjacent.

$\bigskip $

$\qquad \qquad \qquad \ k-2\qquad k-1\qquad k\qquad k+1\qquad k+2$

$x^{2}=(1...1\qquad \ \ 1$ $\ \ \qquad \ \ 1$ $\ \ \ \qquad 0\ \qquad \ \ 1$ 
$\ \ \ \qquad \ \ 0\qquad \qquad 00...0)$

$x^{3}=(1...1\qquad \ \ 1$ $\ \ \qquad \ \ 1$ $\ \ \ \qquad 0$ $\qquad \ \
0\ \ \ \ \qquad \ \ 1\qquad \qquad 00...0)$

$x^{4}=(1...1\qquad \ \ 0$ $\ \ \qquad \ \ 1$ $\ \ \ \qquad 1\ \qquad \ \ 1$ 
$\ \ \ \qquad \ \ 0\qquad \qquad 00...0)$

\bigskip

Since $k-1$ of the first $k$ digits of $x^{2},...,x^{p}$\ must be $1$, and
no two of these \ vertices have the same first $k$ digits, $p$ must satisfy $%
p\leq k+1$. But $p=d-k+2$ implies that $d-k+2\leq k+1.$ A little algebra
gives $\frac{d+1}{2}\leq k$. However, $k\leq \frac{d}{2},$ which implies $%
\frac{d+1}{2}\leq \frac{d}{2}$ a contradiction. \ \ \ $\blacksquare $

\bigskip

\bigskip

\textbf{Proposition 6} \textit{For} $1\leq k\leq \frac{d}{2}:$

(a) $G_{d,k}$\ \textit{is }$(d-1)$-\textit{connected.}

(b) $G_{d,k}$\textit{\ is hamilton connected.}

\bigskip

\textbf{Proof}. (a)\ Balinski's Theorem [12] tells us that every $d$%
-dimensional polytope is $d$-connected. Since $\Delta _{d,k}$ is a $(d-1)$%
-dimensional polytope, it must be $(d-1)$-connected.

(b) Naddef and Pulleyblank [10] proved that if the graph of a $(0,1)$%
-polytope is bipartite, then it is a hypercube. Moreover, if the graph is
nonbipartite, then it is hamilton connected. Proposition 5 implies that $%
G_{d,k}$\ contains $K_{d-k+1}$ as a subgraph. Since $d-k+1\geq 3$, $G_{d,k}$
contains an odd cycle. Therefore, $G_{d,k}$ is not bipartite, and hence, is
hamilton connected. \ $\blacksquare $

\bigskip

\bigskip

\textbf{Proposition 7} \textit{For} $1\leq k\leq \frac{d}{2}$, $G_{d,k}$\ 
\textit{decomposes into} $G_{d-1,k}\cup G_{d-1,k-1}\cup E$, \textit{where} $%
E $ \textit{is a subgraph containing} $\frac{(d-1)!}{(k-1)!(d-k-1)!}$ 
\textit{edges that link} $G_{d-1,k}$ \textit{to} $G_{d-1,k-1}$.

\bigskip

\textbf{Proof}. Consider the subset of $(x_{1}x_{2}...x_{d})$ $\in
V(G_{d,k}) $ that satisfy $x_{1}=1$. These vertices must all satisfy $%
\dsum\limits_{i=2}^{d}x_{i}$ $=k-1$. Let $H_{1}$ be the subgraph induced by
these $\left( 
\begin{array}{c}
d-1 \\ 
k-1%
\end{array}%
\right) $ vertices. Then $H_{1}$ is isomorphic to $G_{d-1,k-1}$. For given
any vertex $x$ in $H_{1}$ we can remove the first coordinate to obtain a
vertex $x^{\prime }$ in $V(G_{d-1,k-1})$. Moreover if $x$ and $y$ are
adjacent in $G_{d,k}$, then $x\cdot y=k-1$. The corresponding vertices $%
x^{\prime }$and $y^{\prime }$ in $G_{d-1,k-1}$ will be adjacent in $%
G_{d-1,k-1}$ since $x^{\prime }\cdot y^{\prime }=k-2.$

\bigskip

Now consider the subset of $V(G_{d,k})$ that satisfy $x_{1}=0$. These
vertices must all satisfy $\dsum\limits_{i=2}^{d}x_{i}$ $=k$, so there are $%
\left( 
\begin{array}{c}
d-1 \\ 
k%
\end{array}%
\right) $ such vertices. Let $H_{0}$ be the subgraph induced by these
vertices. Then an argument similar to the above shows that $H_{0}$ is
isomorphic to $G_{d-1,k}$.

The formula for the number of edges in $G_{d,k}$ given in Proposition 3 can
be used to obtain the equation below, which can then be used to find $%
\left\vert E\right\vert .$

\[
\frac{d!}{2(k-1)!(d-k-1)!}=\frac{(d-1)!}{2(k-1)!(d-k-2)!}+\frac{(d-1)!}{%
2(k-2)!(d-k-1)!}+\left\vert E\right\vert 
\]%
$\blacksquare $

\bigskip

\section{Random walks and the expansion of G$_{d,k}$}

\bigskip

We have demonstrated that $G_{d,k}$ is a tractable graph and many of the
well known graph attributes and parameters of the complete graph $K_{d}$ may
be extended to $G_{d,k}$. In [4], [5] and [8] random walks on the graphs of $%
(0,1)$-polytopes were investigated as a potential algorithm for random
generation of combinatorial objects. In the case of the hypersimplex, the
vertices of $G_{d,k}$\ can be used to represent subsets of $\{1,2,...,d\}$
of size $k$ as follows. Given a vertex $x$, $i$ is in subset $S$ if and only
if $x_{i}=1$. The adjacency criterion given in Proposition 1 allows us to
generate a random neighbor. For given a vertex $x$, generate two random
integers between $1$ and $d$, say $r$ and $s$, until $x_{r}+x_{s}=1$. Then
whichever of $x_{r}$ or $x_{s}$ is equal to $1$ we change to $0$, and
whichever is $0$ we change to 1. Starting with any vertex we may repeat this
process a large number of times. The result is a randomly generated vertex
corresponding to subset of size $k$. We note that there are other known
algorithms to generate random subsets of size $k$ (e.g., see [9]) but
advantages of the above algorithm is that it is easy to code and also an
interesting application of a random walk.

\bigskip

Surprisingly perhaps, the success of the above algorithm is known to depend
on the "edge expansion" properties of $G_{d,k}$. Given a graph $G=(V,E)$,
the \textit{edge expansion} of $G$, denoted $\chi (G),$ is defined as 
\[
\chi (G)=min\left\{ \frac{\left\vert \delta (U)\right\vert }{\left\vert
U\right\vert }:U\subset V,\text{ }U\neq \emptyset ,\text{ }\left\vert
U\right\vert \leq \frac{\left\vert V\right\vert }{2}\right\} 
\]%
where $\delta (U)$ is the set of all edges with one end node in $U$ and the
other one in $V-U$. The edge expansion rate for graphs of polytopes with $%
(0,1)$-coordinates has been recently studied and is an important parameter
for a variety of reasons [4]. In the case of random walks on graphs, "good"
edge expansion implies that the process converges to its limiting
distribution as rapidly as possible [4]. It is known that the hypercube, $%
Q_{n}$ has edge expansion $1$, and has been conjectured that all $(0,1)$%
-polytopes have edge expansion at least $1$ [8]. In [5] it was shown that
the the graph $G_{d,k}$\ has expansion rate at least 1.

\bigskip

When a graph is regular, algebraic graph theory [2] can be used to help
study its expansion rate. If $A$ is the adjacency matrix \ of an $n$-vertex
graph $G$, then $A$ has $n$ real eigenvalues which we denote by $\lambda
_{0}\geq \lambda _{1}\geq \cdots \geq \lambda _{n-1}.$ If $G$\ is a regular
graph with degree $r$, then it is known that $\lambda _{0}=r$, and a result
of Cheeger tells us that $\frac{r-\lambda _{1}}{2}\leq \chi (G)\leq \sqrt{%
2r(r-\lambda _{1})}$ (for a proof, see [4])$.$ For example, the eigenvalues
of the adjacency matrix of $K_{d}$ are $d-1,-1,-1,\ldots \ ,\ -1,$ and
hence, $\frac{d}{2}\leq \chi (K_{d})$. By Proposition 2, we know that the
adjacency matrix associated with $G_{d,k}$ has $\lambda _{0}=r=k(d-k)$. To
investigate the expansion rate of $G_{d,k}$ we need the following
proposition [1].

\bigskip

\textbf{Proposition 8} \textit{For} $1\leq k\leq \frac{d}{2},$ \textit{the
eigenvalues of }$G_{d,k}$\ \textit{are given by }$\lambda
_{j}=(k-j)(d-k-j)-j $, \textit{for }$j=0,1,...,k$, \textit{with
multiplicities} $m_{j}=$ $\left( 
\begin{array}{c}
d \\ 
j%
\end{array}%
\right) -\left( 
\begin{array}{c}
d \\ 
j-1%
\end{array}%
\right) .$

\bigskip

\textbf{Proposition 9} \textit{For }$1\leq k\leq \frac{d}{2},$ \textit{the
edge expansion of }$G_{d,k}$\ \textit{satisfies} $\frac{d}{2}\leq \chi
(G_{d,k})\leq \sqrt{2dk(d-k)}.$

\bigskip

\textbf{Proof}. By Proposition 8, we see that $\lambda _{1}$ $%
=(k-1)(d-k-1)-1 $. Since $r=k(d-k)$, we have that $r-\lambda _{1}\ =d$. If
we now apply Cheeger's Theorem, then $\frac{r-\lambda _{1}}{2}=$ $\frac{d}{2}%
\leq \chi (G_{d,2})\leq \sqrt{2r(r-\lambda _{1})}=$\ $\sqrt{2dk(d-k)}$. $\ \
\blacksquare $

\bigskip

This again extends a property of $K_{d}$, and it is interesting to note that
the lower bound $\frac{d}{2}\leq \chi (G_{d,k})$ is independent of $k$. It
was shown in [5] that $1\leq \chi (G_{d,k}),$ which confirms the conjecture
of Mihail for this special case of hypersimplices. Proposition 9 provides an
improved lower bound and that implies the family of graphs $G_{d,k}$\ has
very good expansion. Consequently, the algorithm mentioned above should be
able to efficiently generate good random subsets.

\end{document}